# Comment: Causal Inference in the Medical Area

**Edward L. Korn**

It is an honor to be a discussant to the Morris Hansen Lecture, and a pleasure to be discussing Don Rubin's talk. Dr. Rubin has clarified over the years many of the deep issues relating to causal inference.

Let me start with a story. About 20 years ago when I was teaching at UCLA, I was eating breakfast one morning at my kitchen table, and my two-and-a-half-year-old daughter was in the next room, lying on her back and kicking the wall with her feet. I told her to stop, which she did for a few seconds, and then began again. I told her to stop again, and that I really meant it. The kicking stopped for a longer period this time, maybe 30 seconds, and then started up again. Just then the Whittier–Narrows earthquake hit, 5.9 on the Richter scale. Our 50-year-old house started shaking like crazy. As I was running into the next room to get my daughter, I ran into her running into the kitchen screaming "I'm sorry, Daddy, I'm sorry. I didn't mean to do it!" Which brings me to my first point: causal inference can be tricky.

Causal inference can be tricky not just for small children, but for epidemiologists and biostatisticians, too. As an example, consider hormone-replacement therapy for postmenopausal women. Dozens of observational studies (including case-control studies and cohort studies) had suggested a 40–50% reduction in coronary heart disease (Stampfer and Colditz, 1991). However, the recently reported results of the Women's Health Initiative trial demonstrated that the treatment had an elevated incidence of coronary heart disease (Manson et al., 2003). Now the statisticians who worked on these epidemiologic studies thought they were making a valid causal inference.

*Edward L. Korn is Mathematical Statistician, Biometric Research Branch, National Cancer Institute, Bethesda, Maryland 90824, USA e-mail: korne@ctep.nci.nih.gov.*

In fact, many women took estrogen replacement therapy partly because they believed that it would offer cardiovascular benefits. However, as the large randomized trial demonstrated, this causal inference from the observational data was completely wrong.

Because of the difficulty of doing randomized clinical trials of certain interventions, and the public health importance of whether these interventions work, I have put the ability to perform causal inference on epidemiologic data on the top of my personal list of "practical importance" of causal inference methods (Figure 1). The hormone-replacement therapy example is, of course, not the only example of medical studies where incorrect causal inferences were made. Let me just mention one other: There were many observational studies that suggested beta carotene would reduce lung cancer incidence; see International Agency for Research on Cancer (1998, pages 64–103) for a summary. However, randomized trials of beta carotene supplements showed that it actually increased the risk of lung cancer. In fact, the epidemiologic data were so strong that when the results of the first trial came out (Alpha-Tocopherol Beta Carotene Cancer Prevention Study Group, 1994), an editorial suggested the possibility that trial results might be due to an "extreme play of chance" (Hennekens, Buring and Peto, 1994). However, after the results of the second trial also showed beta carotene was causing an increase in lung cancer (Omenn et al., 1996), it became clear that the epidemiologic studies had been wrong. To the extent that Dr. Rubin's work can lead to better causal inferences with epidemiologic data of these sorts, it would be of tremendous practical importance.

A cynical colleague of mine suggested that one should not give a discussion like this without mentioning some of your own work. So as an aside, I want to briefly mention a causal analysis I did a few years ago. We were interested in estimating the effect of an orthodontic treatment from observational data (Figure 2). These data were from the University of the Pacific orthodontic clinic, so which orthodontist saw which patients could be assumed to be random. What was definitely not random was which patients received the extraction treatment and which received the nonextraction treatment, because this decision depends on the patient characteristics. Because the treatment decision is not random, one cannot just compare the outcomes for patients who received extraction with those who did not. Instead, we asked each orthodontist to evaluate the other orthodontist's patients and to state which treatment he would have used for that patient. (These evaluations were based only on pretreatment patient records.) The stated orthodontist preferences enabled us to stratify the patient population based on orthodontist preferences, and perform an analysis restricted to the patient subsets where the orthodontists disagreed (strata 2 and 3 in Figure 3). This yielded an appropriate causal inference (Korn, Teeter and Baumrind, 2001). The strata here are similar to Dr. Rubin's principal strata, although we were in a simpler situation because we did not have to estimate which individuals were in which strata using latent variables, but could just observe them. Therefore, we did not have to make the distributional assumptions that seem to be required by Dr. Rubin to make sharp causal inferences.

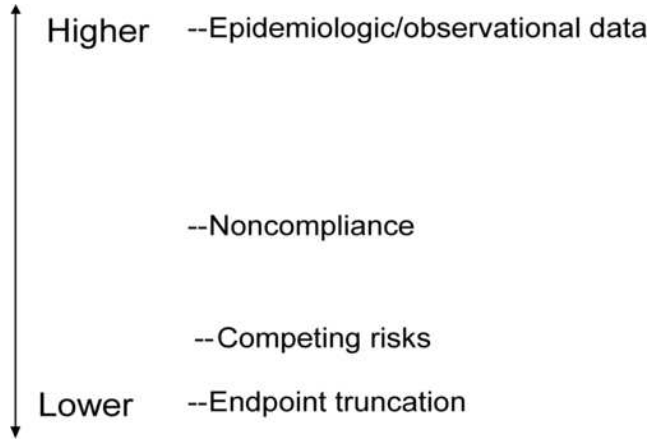

Fig. 1. *Practical importance of causal inference issues in the medical area (personal view).*

Returning to our main discussion and moving down the scale of practical importance, we have noncompliance in randomized clinical trials (Figure 1): Assume that you are doing a randomized trial of a new agent for cancer versus a standard treatment. Now, not all the patients randomized to the new agent may actually take it—some may not take it because they are too sick, some may not take it because they are having bad side effects, and there could be other reasons. An analysis of *effectiveness* analyzes the results of all the patients, based on the treatment arms they were randomized to. This is sometimes called an "intent-to-treat" analysis. A causal-type analysis might be interested in what is sometimes called *efficacy*, the treatment difference that would have been observed if there had been no noncompliance.

The usual arguments for using effectiveness are it is straightforward to estimate with no assumptions, and it estimates in the real world how well the treatment is going to work. The usual arguments for using efficacy are that it estimates the biological effect of treatment better than effectiveness, and it may estimate future effectiveness better than current effectiveness does. This last argument is that if trials show that a new treatment works, then in the future patients may be more compliant with that treatment. This argument, however, leads to a potential issue with Dr. Rubin's methods. His method seems to be estimating the treatment effect only on the subset of patients who would comply given either treatment *in this trial*—not in some future setting where the results of this trial are known. Therefore, one of the usual reasons for being interested in efficacy and not effectiveness seems to be negated by the proposed analysis.

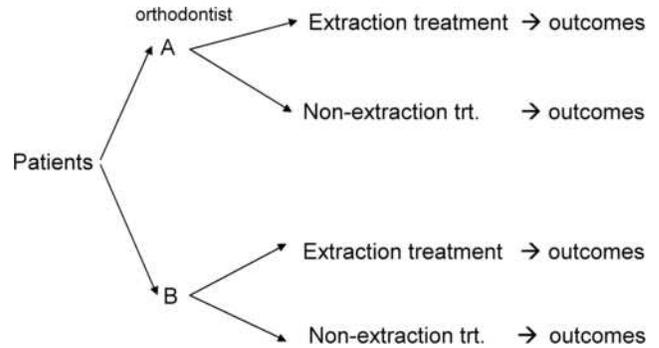

Fig. 2. *Observational orthodontic data.*

|  | orthodontist | |
|---|---|---|
| Patient | A | B |
| Strata | prefers | |
| 1 | X | X |
| 2 | X | N ← |
| 3 | N | X ← |
| 4 | N | N |

Fig. 3. *Observational orthodontic data showing stratification by orthodontist preferences (arrows designate strata used for causal analysis).*

Moving down the scale of practical importance, we come to endpoint truncation (Figure 1). (I should note that words like "truncation" and "censoring" have technical meanings to biostatisticians, and I am not using the word "truncation" here in its technical sense. Instead, I am using it to designate what Dr. Rubin was referring to in his quality-of-life example.) Imagine we are conducting a randomized trial for treating lung cancer, and the outcome is patient-assessed quality of life at 12 months after randomization. Some patients are unfortunately going to die of lung cancer before reaching 12 months, so how does one account in the analysis for a patient who dies at 6 months? You could omit the patient from the analysis, but this leads to obvious bias. You could try to estimate what the patient's quality of life would have been at 12 months if he had not died. But this sounds pretty meaningless, almost supernatural. You could restrict the inference to the subset of patients who would be alive at 12 months regardless of which treatment they were given. This is Dr. Rubin's approach. Finally, you could assign the patient a quality-of-life score consistent with being dead. For example, if you were doing a rank-sum test comparing the 12-month quality-of-life scores between the treatment arms, you could give individuals who died before 12 months the lowest possible score, such as 0. As mentioned by Dr. Rubin, there are scaling issues here, but there are always scaling issues with quality-of-life data, with many being more difficult than this. Because I like this last approach (of assigning a low score to individuals who have died), which does not involve any causal issues, I have put endpoint truncation low on my list of causal issues of practical importance (Figure 1).

There is another kind of endpoint truncation where a causal method might be of more practical importance. Suppose you are conducting a randomized trial of various types of local radiation for head and neck cancer. (Local radiation means radiation just at the tumor site.) As a secondary analysis, one might be interested in how the different types of radiation affect local control of the tumor. One might use a survival analysis with the endpoint being the amount of time from randomization to local recurrence of the tumor. In this type of survival analysis, there are standard ways to handle (i) individuals who die from causes unrelated to their cancer and (ii) individuals who are alive with no evidence of cancer when the analysis is performed. What is more difficult is how to accommodate individuals who have a metastatic development of their cancer, and possibly die, without ever having a local recurrence. This is known as a competing risks problem. I thought about how I would apply Dr. Rubin's methods to this problem. As a first pass, it seems like one would conceptually restrict attention to patients who would have a local recurrence if they were given either randomized treatment, but since there is time involved, it was not obvious how to think about it. Perhaps Dr. Rubin has explored this in some of his other papers.

In summary, I think Dr. Rubin's methods may be very important for interpreting epidemiologic data. I would like to see him take some of the old epidemiologic studies that we now know came to wrong conclusions, apply his methods, and show his methods lead to the right conclusions. That would be impressive.